\documentclass{trbunofficial}
\usepackage{graphicx}
\usepackage{enumitem}    
\usepackage{amsfonts}    
\usepackage[hidelinks]{hyperref}
\usepackage{eurosym}

\usepackage{graphicx}
\usepackage{subcaption}
\usepackage{float}
\usepackage{booktabs}
\usepackage{makecell}

\usepackage{xcolor} 
\definecolor{steelblue}{RGB}{70,130,180} 

\usepackage{amsthm}
\newtheorem{proposition}{Proposition}

\AuthorHeaders{Hammerl, Jin, Seshadri, Rasmussen and Nielsen}
\title{Risk Aware Reservoir Control for Safer Urban Traffic Networks}

\author{%
  \textbf{Alexander Hammerl, Corresponding Author}\\
  asiha@dtu.dk\\
  \hfill\break%
  \textbf{Wen-long Jin}\\
  wjin@uci.edu
  \hfill\break%
  \textbf{Ravi Seshadri}\\
  ravse@dtu.dk
  \hfill\break%
  \textbf{Thomas Kjær Rasmussen}\\
  tkra@dtu.dk
  \hfill\break%
  \textbf{Otto Anker Nielsen}\\
  oani@dtu.dk
}

\newwrite\trbwc
\immediate\openout\trbwc=\jobname-words 
\immediate\write\trbwc{8031}            
\immediate\closeout\trbwc

\begin{document}
\maketitle

\section{Abstract}

We present a risk-aware perimeter-style controller that couples safety and efficiency targets in large, heterogeneous urban traffic networks. The network is compressed into two interacting “reservoirs’’ whose dynamics follow the Generalized Bathtub Model, while accidents are described by a self-exciting (Hawkes) counting process whose intensity depends on vehicle exposure, speed dispersion between reservoirs and accident clustering. Accident occurrences feed back into operations through an analytically simple degradation factor that lowers speed and discharge capacity in proportion to the live accident load. A receding-horizon policy minimizes a mixed delay–safety objective that includes a variance penalty capturing risk aversion; the resulting open-loop problem is shown to possess a bang-bang optimum whose gates switch only at accident times. This structure enables an event-triggered MPC that only re-optimizes when new accidents occur, reducing on-line computation significantly. Parameters are calibrated using OpenStreetMap data for metropolitan Copenhagen to analyze traffic dynamics during morning peak commuter demand. Monte-Carlo simulations demonstrate delay savings of up to 30 \% and accident reductions of up to 35\% relative to an uncontrolled baseline, with a transparent trade-off governed by a single risk parameter.

\hfill\break%
\noindent\textit{Keywords}: Model Predictive Control, Generalized Bathtub Model, Traffic Safety, Incident Prevention, Self-Exciting Point Process
\newpage

\section{Introduction}

Road crashes remain a pressing public-health problem and a sizeable economic burden. The World Health Organization \cite{WHO2023} estimates that they claim close to 1.3 million lives each year, while injuring between 20 and 50 million more, and societal costs are estimated in the high billions annually for major industrialized countries.  Beyond these aggregates, accidents trigger secondary congestion, lost productivity and long-term welfare impacts that weigh heavily on urban economies.

Urban areas are disproportionately affected.  
Dense road grids, frequent intersections and multimodal interactions magnify the secondary delays that follow an incident; even a minor crash may propagate queue spill-back across several signal cycles, disrupt public transport reliability and trigger follow-up collisions. 
Such systemic knock-on effects make it difficult to isolate the “pure” safety component from background congestion and call for an integrated modeling framework that treats accidents and traffic dynamics in a unified manner.

Traffic flow theory provides the foundation for modeling vehicular dynamics and evaluating the effectiveness of accident response measures. Traditional approaches in traffic flow theory focus on spatiotemporal dynamics along individual road segments. This comprises both macroscopic frameworks that treat traffic as a continuous medium \cite{Lighthill1955,Richards1956,Payne1971,Aw2000, Zhang2002} and microscopic car-following models that capture individual vehicle behavior \cite{Bando1995,Treiber2000,Helbing2001,Kesting2010}.
For network-level accident analysis, reservoir-based traffic flow models offer a complementary perspective by describing the dynamics of spatially aggregated traffic quantities across entire networks \cite{Geroliminis2008,Daganzo2008,Helbing2009, Hammerl2024a, Hammerl2025}. These models typically assume a time-invariant relationship between vehicle accumulation and average network flow, known as the Network Fundamental Diagram (NFD). Advanced reservoir frameworks also incorporate relationships between network accumulation and trip completion rates \cite{Vickrey1991,Daganzo2007,Arnott2013,Fosgerau2015}, while the Generalized Bathtub Model \cite{Jin2020,Jin2021,Jin2025} further extends this approach by accounting for the distribution of remaining trip lengths among individual vehicles.
This spatial extension proves particularly relevant for understanding distributional effects that influence effective accident mitigation strategies in complex urban networks. 

On the other end of the spectrum, research on accident modeling and prevention spans a wide spectrum of temporal and spatial aggregation, and the relationship between traffic operations and safety has been a longstanding concern in transportation research. At the facility level, Safety Performance Functions (SPFs) have been widely used to predict crash frequencies based on long-term aggregates of daily traffic volumes and roadway features \cite{AASHTO2010, Gooch2016, Donnell2014, Donnell2016}. While such models typically report a positive correlation between traffic volumes and crash counts, this association may not persist under real-time or spatially disaggregated conditions. For example, positive correlations have also been found between average density indicators (such as the volume-to-capacity ratio) and crash frequencies \cite{Ivan2000,Lord2005,Zhou1997}, while findings regarding the role of average speed remain mixed \cite{Taylor2000, Taylor2002,Elvik2004,Nilsson2004,Aarts2006,Quddus2013}.

Another stream of safety literature investigates the influence of geometric design and intersection typologies. Roundabout studies such as \cite{Ausroads1993} and \cite{KimChoi2013} have demonstrated substantial safety gains, often independent of speed effects. \cite{Qin2011} and \cite{Turner2006} further developed accident prediction models highlighting the interplay between geometric complexity, vehicle familiarity, and traffic flow. At signalized intersections, dilemma-zone behavior \cite{Papaioannou2007} introduces unique risks, which are expected to be mitigated by the rule-based behavior of Connected and Autonomous Vehicles (CAVs).

Microsimulation-based safety assessments have enabled richer representations of conflict interactions. \cite{Carbaugh1998}, \cite{Tibljas2018}, and \cite{Rahman2019} used VISSIM and the Surrogate Safety Assessment Module (SSAM) to evaluate the trajectory-based safety performance of CAVs. These studies confirm that CAV operation reduces surrogate risk measures under various conditions, especially in platooning and mixed-fleet scenarios. \cite{Mahajan2020} proposed a data-driven framework to estimate lane-change-related crash risk based on naturalistic driving data, while \cite{Langer2023} demonstrated the power of high-fidelity crash simulations calibrated with real-world data for urban crash occurrence modeling. These are complemented by qualitative analyses such as \cite{Kinando2018}, which emphasize CAVs' potential to prevent intersection-related conflicts.

Machine learning-based approaches have also gained traction in real-time crash prediction. \cite{Theofilatos2019} compared shallow and deep learning architectures to forecast crash likelihood from traffic sensor data, illustrating the potential of data-driven models in operational safety monitoring. \cite{Dubey2025} proposed a multi-agent framework (FAIRLANE) for adaptive lane prioritization based on safety metrics in mixed autonomy traffic. Relatedly, \cite{Chen2010} reviewed broadcasting strategies for vehicular networks, highlighting how timely dissemination of hazard information can help mitigate collision risk.

At the network level, the Macroscopic Safety Diagram (MSD) framework introduced by \cite{Alsalhi2018} integrates operational state variables (e.g., mean network density) with crash propensity. These studies reveal that crash risk tends to peak in congested conditions—often at densities above the critical value that maximizes flow. Theoretical models based on the two-fluid theory have been proposed to support such observations \cite{Dixit2011,Dixit2013}, although their applicability to urban network dynamics remains limited.

Our work extends the existing literature by developing a network-level control framework that addresses accident risk in metropolitan areas holistically while remaining practical for real-world deployment. The framework requires only aggregate flow measurements and incident reports—data routinely available to traffic management agencies—and captures feedback from live accidents through an analytically simple capacity degradation term. This yields closed-form optimal policies with a bang-bang structure that switches only at accident events, enabling both conceptual transparency and computational efficiency suitable for real-time implementation.

The rest of the paper is structured as follows. Section 2 presents the model description, introducing the Generalized Bathtub Model for traffic dynamics and the self-exciting Hawkes process for accident modeling, along with the capacity degradation mechanism that couples safety and operational performance. Section 3 formulates the risk-aware control problem as a receding-horizon optimization with a mean-variance objective that balances delay minimization against accident risk. Section 4 derives key analytical results, including the well-posedness of the closed-loop system, steady-state optimality conditions, and the bang-bang structure of optimal controls that switch only at accident events. Section 5 demonstrates the framework's practical applicability through numerical experiments on a realistic case study of metropolitan Copenhagen, showing significant improvements in both efficiency and safety metrics. Section 6 concludes with a summary of contributions and discusses promising directions for future research.

\section{Model Description}
\label{sec:framework}

The generalized bathtub model \cite{Jin2020,Jin2021} governs the dynamics of active trips $\delta(t,x)$ in a traffic network with homogeneous conditions, where $\delta(t,x)$ represents the number of trips with remaining distance of at least $x$ at time $t$. The average vehicle speed 
\begin{equation} \label{eq:speed} 
v(t) = V(\rho(t)) 
\end{equation}
only depends on the vehicle density $\rho(t) = \delta(t)/L$, where $\delta(t) = \delta(t,0)$ denotes the total number of active trips and $L$ represents the total lane length in the network's reservoir representation. Without loss of generality, we normalize $L$ to unity in this analysis. The macroscopic speed-density relationship $V(\rho)$ constitutes a time-independent, non-increasing, and non-negative function of density $\rho$.

We define the network fundamental diagram (NFD) through the flow-density relationship
\begin{equation} \label{eq:nfd} 
Q(\rho) = \rho \cdot V(\rho), 
\end{equation}
where $Q(\rho)$ represents the network flow. This relationship is concave $\rho$, attaining its maximum value $q_{\max}$ at the critical density $\rho_{\text{crit}}$, such that $Q(\rho_{\text{crit}}) = q_{\max}$.

Common formulations for the speed-density relationship include:
\begin{align}
V(\rho) &= \begin{cases} 
v_f, & \text{if } \rho \leq \rho_c, \\[0.5ex]
\displaystyle\frac{w(\rho_j - \rho)}{\rho}, & \text{if } \rho > \rho_c
\end{cases} 
&& \text{(triangular MFD, highways)} \\[2ex]
V(\rho) &= \begin{cases} 
v_f, & \text{if } \rho \leq \rho_1, \\[0.5ex]
\displaystyle\frac{q_{\max}}{\rho}, & \text{if } \rho_1 < \rho < \rho_2, \\[0.5ex]
\displaystyle\frac{w(\rho_j - \rho)}{\rho}, & \text{if } \rho \geq \rho_2
\end{cases} 
&& \text{(trapezoidal MFD, urban traffic)}
\end{align}

The parameters $v_f$, $w$, $\rho_c$ (or $\rho_1$), $\rho_j$, and $q_{\max}$ denote the free-flow speed, congestion wave speed, critical density at maximum flow, jam density where speed vanishes, and maximum achievable flow, respectively.

Let $f(t,x)$ represent the exogenous rate at which trip with length of at least $x$ are loaded into the network at time $t$, and let the aggregate exit flow (trip completions) from reservoir $r$ be denoted by $g_r(t)$. The conservation of trips yields the discrete-time formulation
\begin{equation} \label{eq:cl_discrete} 
\delta(t + \Delta t, x) = \delta(t, x + v(t)\Delta t) + f(t,x)\Delta t 
\end{equation}
for infinitesimal time intervals $\Delta t$. Taking the limit $\Delta t \to 0$ and substituting equation \eqref{eq:speed} produces the governing partial differential equation
\begin{equation} \label{eq:cl_cont} 
\frac{\partial \delta}{\partial t}(t,x) - V(\delta(t)) \frac{\partial \delta}{\partial x}(t,x) = f(t,x).
\end{equation}

We partition the metropolitan area under consideration into two interacting reservoirs $\mathcal{R} = \{A,B\}$ with cumulative lane lengths $L_r$ for $r \in \mathcal{R}$. The state variable $N_r(t)$ represents the number of active trips currently located in reservoir $r$ and evolves continuously in time $t \in [0,T]$, where $T$ is the planning horizon. Throughout the analysis we normalize lengths by the total lane supply $\sum_r L_r$, which yields dimensionless densities $\rho_r(t) = N_r(t)/L_r$ bounded by $0 \le \rho_r(t) \le \rho_{j,r}$.

Let $D_r(t)$ denote the (unknown) true entry flow into
reservoir $r$, and let $\widehat D_r(t)$ be the value predicted by the
short‑term demand estimator embedded in the MPC.  
We posit a bounded forecast error
\(
   |D_r(t)-\widehat D_r(t)|\le\varepsilon_D,
\)
and we assume that the true demand itself is uniformly bounded,
\begin{equation}
  0 \;\le\; D_r(t) \;\le\; \overline D_r,
  \qquad \forall\,t\in[0,T],
\end{equation}
where the constant $\overline D_r$ is a known
upper envelope derived, for example, from historical peak counts. No other restrictions are imposed at this stage; the specific forecasting methodology 
will be discussed in the MPC implementation section.

\noindent Let $\varOmega=[\omega_{rs}(t)]_{r,s\in\mathcal R}$ be an exogenous origin-destination (OD) share matrix, where $\textstyle\sum_{s}\omega_{rs}(t)=1$. For intra-district trips ($A\!\to\!A$ and $B\!\to\!B$) we allow a fraction to detour through the other reservoir when that route is faster. We define the detour fractions as $\delta_{A}(t) = \frac{1}{1+\exp\!\bigl[\kappa_{\delta}\,(v_A(t)-v_B(t))\bigr]}$ and $\delta_{B}(t) = \frac{1}{1+\exp\!\bigl[\kappa_{\delta}\,(v_B(t)-v_A(t))\bigr]}$, with elasticity $\kappa_{\delta}>0$. Here $\delta_{A}(t)$ is the share of OD-$A\!A$ demand that chooses the path $A\!\to\!B\!\to\!A$ at time $t$, and analogously for $\delta_{B}(t)$.

The instantaneous flow request from $A$ to $B$ is then given by $d_{AB}(t)= \bigl(\omega_{AB} + \delta_{A}\,\omega_{AA}\bigr)\,D_A(t)$, while $d_{BA}(t)= \bigl(\omega_{BA} + \delta_{B}\,\omega_{BB}\bigr)\,D_B(t)$. This formulation captures how effective inter-district demand incorporates both direct origin-destination flows and detour flows from intra-district trips that choose to route through the alternative district when conditions are favorable. Setting $\kappa_{\delta}=0$ (or $\delta_{r}\equiv 0$) recovers the no-detour baseline without affecting the subsequent analysis.

The control variables are admissible gated flows $u_{rs}(t)$, constrained by $0 \;\le\; u_{rs}(t) \;\le\; \min\!\bigl\{\;\bar u_{rs},\;d_{rs}(t)\bigr\}$ for $(r,s)\in\{(A,B),(B,A)\}$, where $\bar u_{rs}$ is the physical metering capacity and $d_{rs}(t)$ the effective demand defined above.

The resulting conservation equations are \begin{equation}\label{eq:dNa}\dot N_A = D_A(t)\;-\;g_A(t)\;-\;u_{AB}(t)\;+\;u_{BA}(t)\end{equation} and \begin{equation}\label{eq:dNb}\dot N_B = D_B(t)\;-\;g_B(t)\;+\;u_{AB}(t)\;-\;u_{BA}(t)\end{equation}. At this stage no specific functional form is imposed on $V_r(\rho)$.

To model the occurrence of accidents within this framework, let $N_r^{\text{acc}}(t)\in\mathbb N$ be the cumulative number of reported accidents in reservoir $r$ up to time $t$. Under the assumption of independent increments, the number of accidents in the network over any time interval would follow a Poisson distribution with parameter equal to the integrated intensity, even when the accident rate is inhomogeneous in time. However, empirical evidence \cite{Li2018}, \cite{Kalair2021}, \cite{AlaimoDiLoro2025} demonstrates that this independence assumption fails to capture the clustering behavior observed in real traffic networks, where the occurrence of an accident significantly increases the likelihood of subsequent accidents in its vicinity. To account for this empirically observed clustering behavior, we adopt a self-exciting Poisson process with intensity
\begin{equation}
\label{eq:intensity}
  \lambda_r(t) \;=\;
     \underbrace{\alpha_r N_r(t)}_{\text{exposure}}
     \;+\;
     \underbrace{\beta_r\!\int_{0}^{t}
          e^{-\gamma_r(t-s)}\,\mathrm dN_r^{\text{acc}}(s)}_{\text{clustering}}
     \;+\;
     \underbrace{\eta_r\bigl|v_A(t)-v_B(t)\bigr|}_{\text{speed dispersion}},
\end{equation}
where $\alpha_r,\beta_r,\gamma_r,\eta_r\ge0$ are calibration parameters.

To reflect the operational degradation caused by accidents we first define the “live‑accident load’’
\[
   a_r(t)=\sum_{t_i<t} e^{-\gamma_r(t-t_i)},
\]
where the sum runs over the occurrence times \(t_i\) of past accidents in reservoir \(r\). The dimensionless capacity reduction factor is then
\begin{equation}
\label{eq:degradation}
   \chi_r(t)=\frac{1}{1+\kappa_{r}\,a_r(t)},\qquad
   0<\chi_r(t)\le 1,
\end{equation}
with \(\kappa_{r}\ge0\) governing the sensitivity of traffic flow to the prevailing accident load. A larger accident intensity \(\lambda_r\) therefore reduces the effective speed and, through the bathtub dynamics, the discharge capacity of reservoir \(r\). The accident-adjusted speed-density relationship becomes
\begin{equation}
V_r(t)=\chi_r(t)\,V_0\!\bigl(\rho_r(t)\bigr),\qquad r\in\{A,B\},
\tag{\ref{eq:speed}$^\prime$}
\end{equation}
where $V_0\!\bigl(\rho_r(t)\bigr)$ denotes the baseline speed-density function under accident-free conditions.

The clustering term contains a stochastic Itô integral which sums the weighted contributions of all past accidents: each accident at time $s$ adds $\beta_r e^{-\gamma_r(t-s)}$ to the current intensity, with exponential decay capturing how the influence of past accidents diminishes over time. From the definition of these variables, the Kolmogorov Forward System follows immediately for all $N_r^{\text{acc}} \in \mathbb{N}_0$:
\begin{equation}
\label{eq:kolmogorov}
\frac{\partial p_{n}(t)}{\partial t} = 
\begin{cases}
-\lambda_r(t) p_0(t), & n = 0 \\
\lambda_r(t) p_{n-1}(t) - \lambda_r(t) p_n(t), & n \geq 1
\end{cases}
\end{equation}
where $p_n(t) = \Pr\{N_r^{\text{acc}}(t) = n\}$ denotes the probability of having exactly $n$ accidents in reservoir $r$ by time $t$. This system describes the evolution of the probability distribution over the accident count, with the intensity $\lambda_r(t)$ governing the rate of transitions from state $n$ to state $n+1$. 
The exponential kernel $\exp(-\gamma_r\tau)$ preserves the Markov property after augmenting the state with the conditional intensity $\lambda_r$. To understand this relationship, recall that the Markov property requires that the future evolution of a stochastic process depends only on its current state, not on its entire history. In standard Poisson processes, this property holds naturally. However, self-exciting processes violate this property because the intensity at time $t$ depends on the entire history of past events through the integral $\int_{0}^{t} e^{-\gamma_r(t-s)}\,\mathrm dN_r^{\text{acc}}(s)$.

The exponential kernel $\exp(-\gamma_r\tau)$, which weights how past accidents influence future accident rates, preserves the Markov property after augmenting the state with the conditional intensity $\lambda_r$. To understand this relationship, recall that the Markov property requires that the future evolution of a stochastic process depends only on its current state, not on its entire history. In standard Poisson processes, this property holds naturally. However, self-exciting processes violate this property because the intensity at time $t$ depends on the entire history of past events through the integral $\int_{0}^{t} e^{-\gamma_r(t-s)}\,\mathrm dN_r^{\text{acc}}(s)$.

The key insight is that the exponential kernel allows us to compress this entire history into a single state variable. The conditional intensity $\lambda_r(t)$ can be updated recursively via the differential equation: $\frac{d\lambda_r}{dt} = -\gamma_r(\lambda_r - \alpha_r N_r - \eta_r|v_A - v_B|) + \beta_r\frac{dN_r^{\text{acc}}}{dt}$. This equation shows that to compute $\lambda_r(t+dt)$, we need only know the current values of $\lambda_r(t)$, $N_r(t)$, and $N_r^{\text{acc}}(t)$ — not the entire accident history. By augmenting the state space with $\lambda_r$, we transform a non-Markovian process with infinite memory into a Markovian process with finite-dimensional state. This property is essential for computational tractability as it enables us to treat the augmented system as a continuous Markov chain, which greatly simplifies the control formulation in Section \ref{sec:control} and the numerical experiments in Section \ref{sec:numerical}.

Collecting all dynamic variables, the continuous-time state vector is
\begin{equation}
  \mathbf x(t)
     \;=\;
     \Bigl(
       N_A,\;N_B,\;
       \lambda_A,\;\lambda_B
     \Bigr)^{\!\top}\!(t)
     \;\in\;\mathbb R_{\ge0}^{4}.
\end{equation}
Each component is continuously differentiable except $N_r^{\text{acc}}$, which evolves by pure jumps and appears only inside the intensity~\eqref{eq:intensity}. No additional states arise from the routing or detour model because the logistic shares depend solely on the speeds already implicit in $\mathbf x$.

The decision variables are the gated inter-reservoir flows $\mathbf u(t)=(u_{AB},u_{BA})^{\top}(t)$. At every instant they must satisfy
\begin{equation}
  0
  \;\le\;
  u_{rs}(t)
  \;\le\;
  \min\bigl\{\bar u_{rs},\;d_{rs}(t)\bigr\},
  \qquad (r,s)\in\{(A,B),(B,A)\},
\end{equation}
where $\bar u_{rs}$ denotes metering capacity and $d_{rs}(t)$ the effective demand. Define the admissible set $\mathcal U(t)=\bigl\{\mathbf u(t):\text{constraints above hold}\bigr\}$.

We measure efficiency through the running delay proxy $C_T(t)=N_A(t)+N_B(t)$ and safety through the cumulative accident count $N^{\text{acc}}(T)=N_A^{\text{acc}}(T)+N_B^{\text{acc}}(T)$. To capture risk aversion intuitively and with analytical convenience, we adopt a mean-variance surrogate:
\begin{equation}
\label{eq:objective}
  J(\mathbf u)
  \;=\;
  \mathbb E\!\Bigl[
     \int_{0}^{T} c_T\,C_T(t)\,\mathrm dt
  \Bigr]
  +\;
  \underbrace{c_S\,
     \mathbb E\!\bigl[N^{\text{acc}}(T)\bigr]}_{\text{expected safety cost}}
  +\;
  \underbrace{\theta\,
     \mathrm{Var}\!\bigl[N^{\text{acc}}(T)\bigr]}_{\text{risk penalty}}\!,
  \qquad\theta\ge0.
\end{equation}
Setting $\theta=0$ recovers a risk-neutral formulation. Because first and second moments of a linear Hawkes process admit closed-form ODEs (see Section~\ref{sec:analytical}), the augmented optimization problem $\min_{\mathbf u(\cdot)\in\mathcal U(\cdot)} J(\mathbf u)$ remains tractable and retains the analytical structure exploited in Section~\ref{sec:analytical}.
Table~\ref{tab:parameters} summarizes the main symbols used throughout this study.

\begin{table}[H]
  \centering
  \caption{Main symbols used in the model}
  \label{tab:parameters}
  \renewcommand{\arraystretch}{1.1}
  \begin{tabular}{@{}lll@{}}
    \toprule
    \textbf{Symbol} & \textbf{Description} & \textbf{Units} \\
    \midrule
    $L_{A},L_{B}$          & Lane length per reservoir          & km        \\
    $D_{r}(t)$             & True demand flow                   & veh/h     \\
    $\widehat D_{r}(t)$    & Forecasted demand                  & veh/h     \\
    $\varepsilon_{D}$      & Max. forecast error                & veh/h     \\
    $\overline D_{r}$      & Absolute demand ceiling            & veh/h     \\
    $N_{r}(t)$             & Active trips (state)               & veh       \\
    $\rho_{r}(t)$          & Vehicle density                    & veh/km    \\
    $v_{f,r}$              & Free-flow speed                    & km/h      \\
    $\rho_{j,r}$           & Jam density                        & veh/km    \\
    $u_{AB},u_{BA}$        & Gated flows (controls)             & veh/h     \\
    $\bar u_{rs}$          & Metering capacity                  & veh/h     \\
    $d_{rs}(t)$            & Inter-district demand request      & veh/h     \\
    $\omega_{rs}$          & Static OD shares (exogenous)       & —         \\
    $\kappa_{\delta}$      & Detour elasticity                  & —         \\
    $\alpha_{r}$           & Accident exposure coeff.           & h$^{-1}$  \\
    $\beta_{r}$            & Self-excitation gain               & —         \\
    $\gamma_{r}$           & Excitation decay                   & h$^{-1}$  \\
    $\eta_{r}$             & Speed-dispersion coeff.            & h$^{-1}$  \\
    $c_{T}$                & Delay cost weight                  & min$^{-1}$\\
    $c_{S}$                & Safety cost weight                 & —         \\
    $\theta$               & Risk-aversion parameter            & —         \\
    $\tau_{c}$             & Control interval                   & s         \\
    $H_{p}$                & Prediction horizon                 & min       \\
    $N_{p}$                & Horizon steps ($H_p/\tau_{c}$)     & —         \\
    $\kappa_{r}$           & Accident impact coefficient           & h \\

    \bottomrule
  \end{tabular}
\end{table}

\section{The Control Problem}
\label{sec:control}
The on-line controller operates in a receding-horizon manner with control interval $\tau_{c}$ (e.g.\ $\tau_{c}=60\,$s). At the $k$-th control instant $t_{k}=k\tau_{c}$ we solve
\begin{equation}
\label{eq:mpc-problem}
  \begin{aligned}
    \min_{\mathbf u(\cdot)}
      \;&
      \mathbb E\!\Bigl[
         \int_{t_{k}}^{t_{k}+H_{p}}
            \!c_T\,C_T(t)\,\mathrm dt
       \;+\;
         c_S\,N^{\text{acc}}(t_{k}+H_{p})
       \;+\;
         \theta\,\mathrm{Var}\!\bigl[N^{\text{acc}}(t_{k}+H_{p})\bigr]
      \Bigr]                                           \\[3pt]
    \text{s.\,t.}\quad
      &\dot{\mathbf x}(t)=f\bigl(\mathbf x(t),\mathbf u(t),D(t)\bigr),
       \quad t\in[t_{k},t_{k}+H_{p}],                                    \\
      &\mathbf u(t)\in\mathcal U(t),                                     \\
      &\mathbf x(t_{k})\;=\;\mathbf x_{k}\ \text{(measured)},            \\
      &\mathbf x(t)\in\mathcal X,\quad
        N_r^{\text{acc}}(t)\in\mathbb N,\quad
        \forall\,t\in[t_{k},t_{k}+H_{p}],
  \end{aligned}
\end{equation}
where $H_{p}=N_{p}\tau_{c}$ is the prediction horizon and $f(\cdot)$ stacks the differential equations \eqref{eq:dNa}–\eqref{eq:dNb}, the intensity \eqref{eq:intensity}, and the moment ODEs used to evaluate the variance term. Only the first-step control $\mathbf u^{\star}(t_{k})$ is applied; at $t_{k+1}$ the horizon slides forward and the optimization is repeated with updated state and demand forecasts. A terminal-cost or terminal-constraint variant can be recovered by replacing the accident terms at $t_{k}+H_{p}$ with a function $V_f(\mathbf x(t_{k}+H_{p}))$ whenever tighter stability guarantees are required. 

The objective \eqref{eq:mpc-problem} contains the variance of the cumulative accident count $N^{\text{acc}}(t) = N_A^{\text{acc}}(t) + N_B^{\text{acc}}(t)$. To keep the optimization problem deterministic, we append two auxiliary state variables that track the first and second moments of this counting process on-line (cf.\ \cite{DaleyVereJones2007}). For this purpose, define $m(t)=\mathbb E\!\bigl[N^{\text{acc}}(t)-N^{\text{acc}}(t_k)\bigr]$ and $s(t)=\mathbb E\!\bigl[(N^{\text{acc}}(t)-N^{\text{acc}}(t_k))^{2}\bigr]$. $Var[N^{\text{acc}}(t)-N^{\text{acc}}(t_k)]$ follows immediately from these definitions. From $\mathrm d\mathbb E[N^{\text{acc}}] / \mathrm dt = \mathbb E[\lambda_A+\lambda_B]$, the ODE expression 
\begin{equation}
  \dot m(t)
    = \lambda_A(t)+\lambda_B(t),
     m(t_k)=0, \label{eq:meanODE}
\end{equation} with initial condition $m(t_k)=0$ follows for $m(t)$. To derive the evolution equation for $s(t)$, we apply Itô's formula for counting processes:
\begin{equation*}
  d\left[(N^{\text{acc}}(t) - N^{\text{acc}}(t_k))^2\right] =
  2(N^{\text{acc}}(t^-) - N^{\text{acc}}(t_k))\, d(N^{\text{acc}}(t) - N^{\text{acc}}(t_k)) +
  \left[d(N^{\text{acc}}(t) - N^{\text{acc}}(t_k))\right]^2.
\end{equation*}
Since counting process increments are binary (either 0 or 1), we have the identity
\begin{equation*}
  \left[d(N^{\text{acc}}(t) - N^{\text{acc}}(t_k))\right]^2 = d(N^{\text{acc}}(t) - N^{\text{acc}}(t_k)),
\end{equation*}
which yields
\begin{equation*}
  d\left[(N^{\text{acc}}(t) - N^{\text{acc}}(t_k))^2\right] =
  \left[2(N^{\text{acc}}(t^-) - N^{\text{acc}}(t_k)) + 1\right] \, d(N^{\text{acc}}(t) - N^{\text{acc}}(t_k)).
\end{equation*}
Taking expectations and noting that
\begin{equation*}
  \mathbb{E}\left[d(N^{\text{acc}}(t) - N^{\text{acc}}(t_k))\right] = [\lambda_A(t) + \lambda_B(t)]\, dt,
\end{equation*}
we obtain
\begin{equation*}
  \frac{d}{dt}\, \mathbb{E}\left[(N^{\text{acc}}(t) - N^{\text{acc}}(t_k))^2\right] =
  \mathbb{E}\left[2(N^{\text{acc}}(t^-) - N^{\text{acc}}(t_k)) + 1\right] \cdot [\lambda_A(t) + \lambda_B(t)].
\end{equation*}
Substituting $m(t) = \mathbb{E}[N^{\text{acc}}(t) - N^{\text{acc}}(t_k)]$ and $s(t) = \mathbb{E}[(N^{\text{acc}}(t) - N^{\text{acc}}(t_k))^2]$, we arrive at
\begin{equation}
\label{eq:secondODE}
  \dot{s}(t) = 2m(t)[\lambda_A(t) + \lambda_B(t)] + [\lambda_A(t) + \lambda_B(t)], \quad s(t_k) = 0.
\end{equation}

Define the extended state vector $\tilde{\mathbf x}(t)=\bigl(\mathbf x(t), m(t), s(t)\bigr)^{\!\top}$. Because $\lambda_r(t)$ is already a component of $\mathbf x(t)$, the ODEs \eqref{eq:meanODE}–\eqref{eq:secondODE} introduce only two additional scalar states, preserving low dimensionality. At the prediction horizon end time $t_{k}+H_p$ we have
\[
  \mathbb E\bigl[N^{\text{acc}}(t_{k}+H_p)\bigr] = m(t_{k}+H_p), 
  \qquad
  \mathrm{Var}\bigl[N^{\text{acc}}(t_{k}+H_p)\bigr]
      = s(t_{k}+H_p) - m^{2}(t_{k}+H_p).
\]
Thus the cost functional in \eqref{eq:mpc-problem} becomes
\begin{equation}
\label{eq:deterministic-cost}
  J_k
  =
  \int_{t_{k}}^{t_{k}+H_p} \! c_T\,C_T(t)\,\mathrm dt
  + c_S\,m(t_{k}+H_p)
  + \theta\,\bigl[\,s(t_{k}+H_p) - m^{2}(t_{k}+H_p)\bigr],
\end{equation}
which depends only on deterministic variables of the augmented system. The MPC problem \eqref{eq:mpc-problem} is therefore equivalent to a finite-dimensional deterministic optimal-control problem with state $\tilde{\mathbf x} \in \mathbb R^{6}_{\ge0}$ and control $\mathbf u \in \mathcal U$. All theoretical properties shown later—existence, bang-bang structure, and performance bounds—carry over unchanged, as the reformulation introduces no additional uncertainty and preserves the affine-in-control Hamiltonian required for the Pontryagin analysis in Section~\ref{sec:analytical}.


The risk-aversion coefficient \( \theta \ge 0 \)
penalizes the variance of accident counts.
Setting \( \theta = 0 \) recovers a risk-neutral policy; larger values trade expected efficiency for a tighter crash distribution. 

Since the risk parameter $\theta$ does not correspond to a directly observable economic quantity, we determine its value through trade-off analysis. We solve problem~\eqref{eq:deterministic-cost} for $\theta \in \{0, 0.05, \ldots, 1\}$ and plot the resulting expected accident count $\mathbb{E}[N^{\text{acc}}]$ versus standard deviation $\mathrm{Std}[N^{\text{acc}}]$. This curve reveals how increasing $\theta$ reduces accident variability at the cost of higher average accidents. We select $\theta$ where the curve exhibits significant curvature—the point beyond which variance reduction becomes disproportionately expensive. For our case study, $\theta = 0.2$ achieves this balance.

To solve the receding-horizon problem in~\eqref{eq:mpc-problem} numerically, we note that with six continuous states and two controls, the problem remains modest in size. We discretize it by direct multiple shooting using a step \( \Delta t = \tau_c \). The resulting nonlinear programme (NLP) contains \( N_p \times 2 = 180 \) control variables and \( N_p \times 6 = 540 \) state continuity constraints for a baseline horizon of \( N_p = 90 \) minutes.
The optimization problem can be solved efficiently using standard nonlinear programming solvers, with computation times well below the control interval requirements. For resource-constrained implementations, the dynamics can be linearized around the current state:
\[
  \tilde{\mathbf{x}}_{j+1}
    = A_j \tilde{\mathbf{x}}_j + B_j \mathbf{u}_j + b_j,
  \qquad
  j = 0, \dots, N_p - 1.
\]
The piecewise affine structure of the exit flow function and the linearity of the moment equations ensure that linearization introduces minimal error. This reformulation yields a convex quadratic program that can be solved extremely rapidly, making the approach suitable for real-time deployment even on embedded systems.

If still faster evaluation is required, one may pre-compute an explicit MPC law by parametric QP, or train a compact neural network on optimal state-control pairs.
Because both surrogates respect the same hard constraints, they inherit the feasibility and stability guarantees proved in Section~\ref{sec:analytical}.

\section{Analytical Results}
\label{sec:analytical}
In this section, we derive analytical properties of the control program that guide the design of an efficient and implementable control policy. Our analysis requires only minimal assumptions about the speed-density relationship $V_r(\rho)$, without specifying a particular functional form such as triangular or trapezoidal relationships. Specifically, we assume that $V_r(\rho)$ is continuously differentiable and strictly decreasing for $\rho > \rho_{c,r}$, that the flow function $Q_r(\rho) = \rho V_r(\rho)$ is concave in $\rho$, and that the standard boundary conditions $V_r(0) = v_{f,r}$ and $V_r(\rho_{j,r}) = 0$ hold. These structural properties suffice to establish the key analytical results that follow.

Before analyzing optimality, we establish that the closed-loop dynamics are well-posed, ensuring that every admissible control generates a unique, bounded trajectory on the finite horizon $[0,T]$. The augmented state vector $\tilde{\mathbf x}(t)=(N_A,N_B,\lambda_A,\lambda_B,m,s)^\top$ evolves as a piecewise-deterministic Markov process. Between accident events, the components $(N_A,N_B,m,s)$ follow the locally Lipschitz ordinary differential equations \eqref{eq:dNa}–\eqref{eq:secondODE}, while the intensity processes satisfy
\[
   \dot\lambda_r(t)= -\gamma_r\lambda_r(t) + \alpha_r\dot N_r(t)
                     + \eta_r\lvert v_A(t)-v_B(t)\rvert ,
   \quad r\in\{A,B\},
\]
which are also Lipschitz on any compact subset of the non-negative orthant.

Accident occurrences follow Poisson processes with time-varying rates $\lambda_r(t)$. To establish boundedness, we note that the intensity satisfies
\[
\lambda_r(t) \le \alpha_r\overline D_r + \eta_r v_{f,\max} + \beta_r t,
\]
where the linear term $\beta_r t$ represents the maximum cumulative self-excitation effect. In the worst-case scenario, an accident could occur in every infinitesimal time interval \( \mathrm{d}t \). By time \( t \), this would amount to at most \( t / \mathrm{d}t \) such events. Since each accident increases the intensity by \( \beta_r \), the total cumulative increase in intensity is bounded above by \( \beta_r \cdot \frac{t}{\mathrm{d}t} \cdot \mathrm{d}t = \beta_r t \). While conservative, this envelope ensures that the intensity remains finite on $[0,T]$, guaranteeing non-explosive counting processes with only finitely many jumps.

The theory of piecewise-deterministic Markov processes \citep{Davis1984} guarantees that the system has a unique solution $\tilde{\mathbf x}(t)$ that evolves continuously between accident events (with jumps only in the intensity at accident times) and remains non-negative and bounded. This mathematical well-posedness translates directly to computational reliability: our numerical solver can find a unique, physically meaningful control trajectory, and the optimization algorithms have well-defined gradients to follow throughout the prediction horizon.

Next, we consider a simplified steady-state regime to provide both a theoretical benchmark and practical initialization for the MPC controller. We consider the simplified setting in which external demands are constant, the risk parameter is set to $\theta=0$, and the horizon is extended to infinity. With $\theta=0$ the objective reduces to minimizing the long-run average cost
\[
  \bar J = c_T\,(N_A^\ast+N_B^\ast)
         + c_S\,\frac{\alpha_A N_A^\ast+\alpha_B N_B^\ast}
                          {\gamma_A+\gamma_B},
\]
where the second term represents the steady-state accident cost obtained by setting $\dot\lambda_r=0$ in \eqref{eq:intensity}. Since both cost components grow linearly with reservoir occupancies, minimizing $\bar J$ is equivalent to minimizing $N_A^\ast+N_B^\ast$.

In \cite{Jin2020}, the relationship
\[
g(t) = \phi(t,0) \cdot n(t) \cdot V(\rho(n))
\]
is derived for the outflow rate \( g(t) \), where \( \phi(t,0) \) denotes the fraction of active trips whose remaining travel distance converges to zero. Furthermore, \cite{Jin2020} demonstrates that \( \phi(t,0) \) is time-invariant if and only if the initial trip length distribution of the inflow follows a negative exponential distribution with scale parameter \( B > 0 \):
\[
f_r(t,x) = \frac{1}{B} \cdot \exp\left(-\frac{x}{B}\right).
\]
Under this condition, \( \phi(t,0) = \frac{1}{B} \) holds, whereby the outflow rate \( g(t) \) can be expressed as a time-independent function determined solely by the instantaneous traffic density:
\[
g_r(N_r) = \frac{1}{B_r} \cdot N_r \cdot v_r(N_r), \qquad r \in \{A,B\}.
\]
In this special case, the bathtub model reduces to the ordinary differential equation
\[
\frac{dn}{dt} = f(t) - \frac{1}{B} \cdot n \cdot V(n).
\]
For the subsequent steady-state analysis, we assume that the distribution of newly entering trips is time-independent and negative exponential, \( f_r(t,x) = F_r e^{-x/B_r} \) with scale parameter \( B_r > 0 \), where \( F_r \) denotes the total inflow to reservoir \( r \). At steady state with $\dot N_A=\dot N_B=0$, the mass-balance relations derived from \eqref{eq:dNa}–\eqref{eq:dNb} yield
\[
  g_A(N_A^{\ast})+u_{AB}^{\ast}-u_{BA}^{\ast}=F_A, \qquad
  g_B(N_B^{\ast})+u_{BA}^{\ast}-u_{AB}^{\ast}=F_B.
\]
Adding these equations shows that the sum of discharges equals the system inflow: $g_A(N_A^{\ast})+g_B(N_B^{\ast})=F_A+F_B$. Thus, while the gating controls can reallocate occupancies between reservoirs, the aggregate outflow remains fixed by demand.

Formulating this as a constrained optimization problem with Lagrangian
\[
  \mathcal L = N_A+N_B -\lambda\,[\,g_A(N_A)+g_B(N_B)-F_A-F_B\,],
\]
the first-order conditions yield the marginal-equality rule
\[
     g_A'\!\bigl(N_A^{\ast}\bigr)
     \;=\;
     g_B'\!\bigl(N_B^{\ast}\bigr).
\]
Since $g_r(N_r)=B_r^{-1}N_r\,v_r(N_r)$ and the speed functions $v_r(\cdot)$ are strictly decreasing, this equation has a unique solution that minimizes total delay.

To implement this optimality condition through the gating controls, we define the net flow offset $\Delta u=u_{AB}^{\ast}-u_{BA}^{\ast}$. From the balance equations, we obtain
\[
  g_A(N_A^{\ast})-g_B(N_B^{\ast}) = F_A-F_B-2\Delta u.
\]
The optimal control strategy directs flow from the reservoir with lower marginal discharge capacity to the one with higher marginal discharge capacity, thereby equalizing the rate at which additional vehicles increase outflow in each reservoir. Within the admissible range $0\le u_{rs}\le\bar u_{rs}$, this yields the bang-bang rule
\[
  u_{AB}^{\ast}= \bigl[\Delta u^{\ast}\bigr]_{+},
  \qquad
  u_{BA}^{\ast}= \bigl[-\Delta u^{\ast}\bigr]_{+},
\]
where
\[
  \Delta u^{\ast}= \frac{F_A-F_B-g_A(N_A^{\ast})+g_B(N_B^{\ast})}{2},
\]
and $[x]_{+} = \max(0,x)$ denotes the positive part operator. This ensures that only non-negative flows are implemented. This closed-form steady-state solution provides two practical benefits. First, it serves as an analytic benchmark: when the prediction horizon is long relative to system time constants, the finite-horizon MPC solution approaches $(u_{AB}^{\ast},u_{BA}^{\ast})$. Second, it supplies a warm-start initial guess that significantly reduces solver runtime, an important consideration for real-time implementation.

Next, we turn to the influence of the risk-aversion parameter \(\theta>0\) that penalizes the variance of the cumulative accident count. To obtain analytic insight we retain the same infinite-horizon, time-invariant setting as above, but augment the
long-run cost by the steady-state variance term
\(\theta\,\mathrm{Var}[N^{\mathrm{acc}}]\).
For a linear Hawkes process with exponential kernel the first two
moments of the count admit closed-form expressions up to lower order terms
\cite{DaleyVereJones2007}:
\[
   \mathbb E[N^{\mathrm{acc}}]
        = \frac{\alpha_A N_A+\alpha_B N_B}{\gamma}
          +\text{const.},
   \qquad
   \mathrm{Var}[N^{\mathrm{acc}}]
        = \frac{\sigma_A N_A+\sigma_B N_B}{\gamma^{2}}
          +\text{l.o.t.},
\]
where \(\gamma>\max\{\beta_A,\beta_B\}\) and
\(\sigma_r:=\alpha_r(\gamma+\beta_r)\).
Hence, up to lower-order terms, the variance is linear in the reservoir occupancies, so the long-run cost can be written as
\[
  \bar J(\theta)
     = \bigl[c_T+\theta\,\sigma_A/\gamma^{2}\bigr]\,N_A
       +\bigl[c_T+\theta\,\sigma_B/\gamma^{2}\bigr]\,N_B
       +\text{const.}
\]
Risk aversion therefore acts as an occupancy surcharge that
inflates the marginal cost of keeping a vehicle inside the system,
with different multipliers for \(A\) and \(B\) whenever
\(\sigma_A\neq\sigma_B\). The steady-state optimization problem remains \(\min_{N_A,N_B\ge0} \bar J(\theta)\)
subject to \(g_A(N_A)+g_B(N_B)=F\).
Forming the Lagrangian
\(
 \mathcal L = \bar J(\theta) - \lambda\,[\,g_A(N_A)+g_B(N_B)-F\,]
\)
and taking first derivatives gives the modified
marginal-equality conditions
\[
   \frac{c_T+\theta\sigma_A/\gamma^{2}}{g_A'(N_A^{\theta})}
   \;=\;
   \frac{c_T+\theta\sigma_B/\gamma^{2}}{g_B'(N_B^{\theta})}
   \;=\;\lambda^{\theta}.
   \tag{ME\(_\theta\)}
\]
Comparing (ME\(_\theta\)) with the risk-neutral rule
\(g_A'(N_A^{0})=g_B'(N_B^{0})\) reveals two qualitative effects:

\begin{itemize}
\item \textit{Uniform bias.}\; Because the risk surcharge inflates the marginal cost of occupancy, the solution to the weighted problem is no longer the minimizer of \(N_A+N_B\).
Consequently the aggregate steady‑state stock increases,
i.e. \(N_A^{\theta}+N_B^{\theta}\ge N_A^{0}+N_B^{0}\).
Whether the extra vehicles are stored predominantly in \(A\) or \(B\) depends on the relative variance coefficients (\(\sigma_A\) versus \(\sigma_B\)); the reservoir with the lower \(\sigma_r\) will absorb the larger share of the additional occupancy.
\item \textit{Directional bias.}\; If \(\sigma_A\neq\sigma_B\) the ratio of surcharges differs, so the controller preferentially relieves the reservoir with the larger variance coefficient, driving the system towards a configuration in which the riskier region carries fewer vehicles.
\end{itemize}

\noindent The gating policy that enforces (ME\(_\theta\)) is again bang-bang.
Let \(\Delta u^{\theta}=u_{AB}^{\theta}-u_{BA}^{\theta}\) and define
the function
\(\Phi(N_A)=g_A'(N_A)/(c_T+\theta\sigma_A/\gamma^{2})\).
Because \(\Phi\) is strictly decreasing, the controller can be
implemented as a two-threshold rule:
increase \(\Delta u\) whenever
\(\Phi(N_A)>\Phi(N_B)\) and decrease it when the inequality is
reversed, subject to the same metering bounds as before.
For small \(\theta\) a first-order expansion gives the approximate
shift
\[
  \Delta u^{\theta}
    = \Delta u^{0}
      - \theta\,\frac{\sigma_A-\sigma_B}{2\gamma^{2}\,
            g_A'(N_A^{0})+2\gamma^{2}\,g_B'(N_B^{0})}
      +\mathcal O(\theta^{2}),
\]
showing explicitly how greater risk sensitivity moves the operating point towards the safer reservoir.  In practice the parameter \(\theta\) thus offers a transparent knob: setting \(\theta=0\) reproduces the purely efficiency-driven solution, while increasing \(\theta\) yields progressively more conservative gating without altering the structural simplicity of the controller.

While the steady-state analysis provides insights into the long-term behavior and the effect of risk aversion, practical implementation requires understanding the transient dynamics of the optimal control. In the MPC framework, at each control instant we solve an open-loop optimization problem—computing the optimal control trajectory over the entire prediction horizon given the current state—but then apply only the first control action before re-solving with updated measurements at the next time step. This receding-horizon approach creates a feedback policy from a sequence of open-loop solutions. We can establish that each of these open-loop solutions possesses a particularly simple structure: each gating input assumes only its minimum or maximum value and changes at most twice along the finite horizon. 

\noindent The argument proceeds via the Pontryagin Maximum Principle (PMP), which transforms the dynamic optimization problem into a two-point boundary value problem by introducing the costates $\boldsymbol\lambda(t)$, which can be interpreted as the marginal value of each state variable along the optimal trajectory. The key insight is that optimal controls must minimize the Hamiltonian $\mathcal{H}(\mathbf{x}, \mathbf{u}, \boldsymbol{\lambda}, t) = L(\mathbf{x}, \mathbf{u}, t) + \boldsymbol{\lambda}^{\top} f(\mathbf{x}, \mathbf{u}, t)$ at each instant in time, where $L$ is the instantaneous cost, $f$ represents the system dynamics, and $\boldsymbol{\lambda}$ are the costate variables. For our system, this becomes:
\[
   \mathcal H = c_T\,C_T(t) + \boldsymbol\lambda^{\top}f\bigl(\tilde{\mathbf x},\mathbf u,t\bigr),
\]
with dynamics $f$ from Section \ref{sec:control}. Because the augmented dynamics are affine (linear plus constant) in the control vector $\mathbf u=(u_{AB},u_{BA})^{\top}$ and the running cost $c_T\,C_T(t)$ does not depend on $\mathbf u$, the Hamiltonian simplifies to a form that is linear in each control component:
\[
   \mathcal H=\boldsymbol\lambda^{\top}f\bigl(\tilde{\mathbf x},t\bigr)
              +\mathbf u^{\top}\Lambda,
\qquad
   \Lambda=
      \begin{pmatrix}
         \lambda_{N_B}-\lambda_{N_A}\\[2pt]
         \lambda_{N_A}-\lambda_{N_B}
      \end{pmatrix},
\]
where $\lambda_{N_r}$ denotes the costate component associated with reservoir occupancy $N_r$.

The switching structure follows directly from PMP because both gating variables enter the Hamiltonian affinely and with opposite signs. Define the switching function \(S_{AB}(t)=\lambda_{N_B}(t)-\lambda_{N_A}(t)\). The optimal input \(u_{AB}^{\star}(t)\) attains its lower bound whenever \(S_{AB}(t)>0\) and its upper bound whenever \(S_{AB}(t)<0\); a sign change of \(S_{AB}\) therefore marks a control switch.  A bound on the number of such sign changes is obtained from the costate equations. According to the PMP framework, the relevant costate dynamics for the augmented deterministic system are
\[
  \dot\lambda_{N_A} = 
     -c_T
     -\bigl(\alpha_A+\theta\sigma_A/\gamma^{2}\bigr)
      \,\lambda_{\lambda_A},
  \quad
  \dot\lambda_{N_B} = 
     -c_T
     -\bigl(\alpha_B+\theta\sigma_B/\gamma^{2}\bigr)
      \,\lambda_{\lambda_B},
\]
where \(\lambda_{\lambda_r}\) is the costate multiplier associated with the accident intensity state \(\lambda_r\). Because \(\lambda_{\lambda_r}\) evolves through the linear ODE \(\dot\lambda_{\lambda_r}=-\gamma_r\lambda_{\lambda_r}\),
it remains continuous and bounded, implying that both
\(\dot\lambda_{N_A}\) and \(\dot\lambda_{N_B}\) are continuous. Subtracting the two equations gives the derivative of the switching function,
\[
   \dot S_{AB}(t)=
      \bigl(\alpha_B+\theta\sigma_B/\gamma^{2}\bigr)
      \lambda_{\lambda_B}(t)
     -\bigl(\alpha_A+\theta\sigma_A/\gamma^{2}\bigr)
      \lambda_{\lambda_A}(t),
\]
which is continuous and—crucially—does not depend on either gating control.  Between any two consecutive state jumps (accident events) the evolution of \(\lambda_{\lambda_r}\) is exponential, so \(\dot S_{AB}\) retains a fixed sign on that interval; consequently \(S_{AB}(t)\) is monotone there and can cross zero at most once.

The control can switch only when the switching function
\(S_{AB}(t)=\lambda_{N_B}(t)-\lambda_{N_A}(t)\) crosses zero. Between two consecutive accident events—i.e.\ on an interval \((\tau_i,\tau_{i+1})\) during which the state is smooth—the costate component \(\lambda_{\lambda_r}\) obeys the linear ODE
\(\dot\lambda_{\lambda_r}=-\gamma_r\lambda_{\lambda_r}\) and therefore remains continuous. As the right‑hand sides in the costate equations depend only on continuous variables, \(\dot S_{AB}(t)\) is continuous
as well, retains a fixed sign on \((\tau_i,\tau_{i+1})\), and \(S_{AB}\)
is consequently monotone there—allowing at most one zero crossing. Each accident starts a new smooth segment and introduces the potential for one further sign change.  Hence, over any finite horizon the number
of switches of \(u_{AB}^{\star}(t)\) is bounded by one plus the number of accident events that actually occur in that horizon. The same
bound holds for \(u_{BA}^{\star}(t)\) because its switching function is \(-S_{AB}(t)\).

The above derivations demonstrate that the optimal online control policy can only change when an accident is registered, thereby justifying the adoption of an \emph{event-triggered} update interval for the MPC. At the initial time, we integrate the adjoint equations once, locate the single zero of \( S_{AB}(t) \), and store the corresponding switch time \( t^{\ast} \). During operation, the gate maintains its current value until an accident is recorded in either reservoir, at which point we recompute a new \( t^{\ast} \) based on the updated state. If no accident occurs, no recomputation is needed, making the approach computationally efficient.

Even in large-scale traffic networks such as Copenhagen's city center (detailed in Section 6), empirical data shows that registered accidents rarely exceed 6 occurrences per hour under typical operating conditions. Given our calibrated Hawkes parameters, the optimization routine is therefore invoked at most 6 times hourly—representing a significant computational burden reduction compared to conventional MPC implementations with 60-second receding horizons, while producing mathematically identical control trajectories. This efficiency gain stems directly from eliminating redundant recomputations during accident-free periods.

\begin{proposition}
Under bounded demand and admissible controls the closed‑loop system is well posed; in steady state the risk‑neutral policy is explicit and the risk‑averse policy follows a weighted marginal balance; over any horizon the optimal gates are bang‑bang with at most
one switch more than the number of accident events. The
event‑triggered MPC therefore achieves near‑optimal delay–safety trade‑offs with modest on‑line computational effort.
\end{proposition}

\section{Numerical Experiments}
\label{sec:numerical}

To demonstrate the practical applicability of the risk-aware MPC framework, we implement the controller on a realistic case study of metropolitan Copenhagen, showing how the bang-bang control structure performs under stochastic demand scenarios and realistic network conditions.

The metropolitan area is divided into two macroscopic reservoirs that capture the morning and evening commuting patterns motivating cordon-style gating control. \textbf{Reservoir A} represents the city core, defined by a polygon that encloses Copenhagen’s main business districts. This region includes the historic city center, Frederiksberg, Vesterbro, Nordhavn, and the northern part of Ørestad. It captures areas characterized by dense office and institutional development, a high concentration of commuters, moderate traffic speeds, and frequent signalized junctions. The area experiences significant pedestrian activity and multimodal interactions, particularly during peak hours, with extensive bus routes and cycling infrastructure competing for limited road space.
\textbf{Reservoir B} comprises the surrounding motorway ring and radial feeders extending approximately 15 km from the city center. This grade-separated network has higher capacity but experiences more severe accident propagation once congestion sets in. The area includes major interchanges, merging zones, and high-speed segments where speed differentials between entering and through traffic exacerbate collision risks. 

Road network data was obtained from OpenStreetMap (OSM) using the OSMnx Python library \cite{Boeing2017}. The administrative boundary of Copenhagen was retrieved via Nominatim geocoder, buffered by 0.01° to include the ring road, and all edges tagged \texttt{highway=drive} were loaded. For each reservoir, we computed: (i) total directed edge length, yielding lane-kilometer supply $L_r$; and (ii) mean posted speed limit, assigned as free-flow speed $v_{f,r}$.

Area A included the full road hierarchy relevant for urban traffic: the primary network (major \texttt{motorway}, \texttt{trunk} and associated links), secondary network (\texttt{primary}, \texttt{secondary} and links), and local distribution (\texttt{tertiary}, \texttt{unclassified}, \texttt{residential}).Area B-A focused exclusively on high-capacity corridors: highway-grade roads only (\texttt{motorway}, \texttt{trunk} and links).

Lane counts were determined using explicit OSM \texttt{lanes} tags when available, otherwise estimated based on Danish infrastructure standards: motorway (4 lanes dual carriageway), trunk roads (4 lanes major arterials), primary roads (3 lanes), secondary/tertiary (2 lanes), and residential (1 lane effective capacity). For computational efficiency, 30\% of urban segments and 50\% of highway segments were randomly sampled, with results scaled proportionally. The query yielded an cumulative lane length of 328.9 km for Network A and 79.5 km for Network B. Maximum flow capacity from the metropolitan region (B-A) to the urban core (A) was estimated using a Fermi-style approach based on the perimeter length of Area A and standard traffic engineering principles. The perimeter of Area A was calculated using geodesic distance measurements, yielding a total boundary length of approximately X.X km.

Entry point density was estimated assuming vehicular access points every 200 meters along the perimeter, resulting in N entry points. Each entry point was assigned an average of 2.5 lanes based on typical Copenhagen street cross-sections at urban boundaries. Theoretical lane capacity was set at 1,800 vehicles per hour per lane, consistent with urban arterial standards.

\subsection{Maximum Perimeter Flow Estimation}

The maximum flow constraint for the flow from region B-A to urban core A was estimated using the perimeter length (22.5 km) and assuming an entry point at every 250m with 2.5 average lanes per point at 1,800 veh/h/lane capacity. This yielded a peak hour capacity of approximately 43,000 veh/h. The perimeter flow is practically unbounded by perimeter capacity because the obtained bound is much higher than demand.

The generalized bathtub model parameters reflect typical characteristics of each network type. Zone A (Urban Core) has free-flow speed of 35 km/h (urban arterials with traffic signals), jam density of 180 veh/km/lane (typical urban congestion), and maximum capacity of 1,580 veh/h/lane (optimistic with good signal coordination). Zone B-A (Metropolitan Ring) has free-flow speed of 85 km/h (highway/trunk roads), jam density of 140 veh/km/lane (highway conditions), and maximum capacity of 2,350 veh/h/lane (optimistic highway capacity). These capacity values are slightly above literature references as they represent idealized conditions without disruptions, with accident impacts separately modeled through the degradation function.

\begin{figure}[H]
\centering
\includegraphics[width=0.8\textwidth]{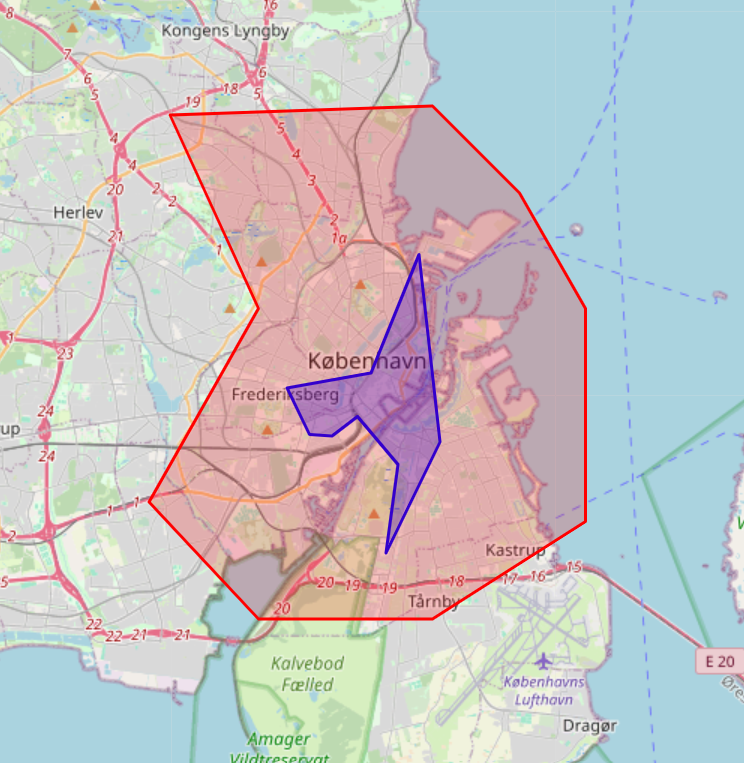}
\caption{Geographic boundaries of the study areas used for model calibration for the urban core area (blue) and the peripheral highway network (red). The area delineation of the city center was conducted considering all essential business districts with pronounced commuter activity to ensure representative coverage of traffic-generating activity centers.}
\end{figure}

The parameterization reflects the empirically established finding that traffic accidents on highways have lower frequency but stronger operational impacts on traffic flow. To analyze the effects of varying accident rates, two parameter sets (low and high disruption intensity) were implemented for the urban core.

\begin{table}[H]
\centering
\renewcommand{\arraystretch}{1.3}
\begin{tabular}{lcccccccccc}
\toprule
\textbf{Reservoir} 
& \multicolumn{2}{c}{$\alpha_r$ (h$^{-1}$)} 
& \multicolumn{2}{c}{$\beta_r$ (–)} 
& \multicolumn{2}{c}{$\gamma_r$ (h$^{-1}$)} 
& \multicolumn{2}{c}{$\eta_r$ (h$^{-1}$)} 
& \multicolumn{2}{c}{$\kappa_r$ (h)} \\
& \textbf{Low} & \textbf{High} 
& \textbf{Low} & \textbf{High} 
& \textbf{Low} & \textbf{High} 
& \textbf{Low} & \textbf{High} 
& \textbf{Low} & \textbf{High} \\
\midrule
City core A 
& 1.5e-4 & 2.25e-4 
& 0.40 & 0.60 
& 1.2 & 1.8 
& 5e-4 & 7.5e-4 
& 0.20 & 0.30 \\
Motorway B 
& 0.5e-4 & 0.5e-4 
& 0.20 & 0.20 
& 1.0 & 1.0 
& 3e-4 & 3e-4 
& 0.35 & 0.35 \\
\bottomrule
\end{tabular}
\caption{Calibration parameters with low and high values. High values for City core A are scaled by 1.5, Motorway B values are unchanged.}
\end{table}

\subsection{Speed-Density Functions for the Triangular Fundamental Diagram}

Based on the given parameters, we can derive the following speed-density functions $V_A(\rho)$ and $V_B(\rho)$ for each zone:

\[
v_A(\rho) =
\begin{cases}
35, & \text{if } \rho \leq 45.14 \\
\frac{11.72(180 - \rho)}{\rho}, & \text{if } \rho > 45.14
\end{cases}
\]
and
\[
v_{BA}(\rho) =
\begin{cases}
85, & \text{if } \rho \leq 27.65 \\
\frac{20.91(140 - \rho)}{\rho}, & \text{if } \rho > 27.65
\end{cases}
\]

The simulated traffic scenario is a morning peak hour with predominantly centripetal traffic flows from the peripheral network to the urban core area. For model simplification, it is assumed that during the morning peak period, no vehicles travel to the peripheral network as their destination. To quantify the traffic shares of origin~A and~B in the total traffic volume, the population within the city center is determined via a WorldPop query as 98{,}838 inhabitants. With a total population of Copenhagen of 667{,}099 inhabitants (2025 census data), this yields an estimated share of the city center in the total traffic demand of 15\%. We assume a constant total inflow of 12{,}000 vehicles per hour over a period of 60 minutes, followed by an abrupt drop to zero. During the subsequent 15-minute clearance phase, no additional vehicles enter the network, but vehicles already in the system can complete their journeys. For modeling the trip length distribution at network entry, we employ the log-normal distribution $\text{LogN}(\mu, \sigma^2)$, which according to~\cite{Martinez2021} provides the best fit to empirically observed trip distance distributions. Parameter estimation for $\mu$ and $\sigma$ is performed via Monte Carlo simulation of the spatial distribution of origin and destination points. This yields the following parameters:

\begin{itemize}
    \item Vehicles entering the network in Reservoir A: \(\mu_A = 2.27~\text{km}, \quad \sigma_A = 1.27~\text{km}\)
    \item Vehicles entering the network in Reservoir B: \(\mu_B = 2.69~\text{km}, \quad \sigma_B = 1.60~\text{km}\)
    \item Vehicles transitioning from Reservoir B to Reservoir A: \(\mu_{B \to A} = 2.79~\text{km}, \quad \sigma_{B \to A} = 1.44~\text{km}\)
\end{itemize}

For the offline optimization of the bang-bang structure of the optimal control derived in Chapter~\ref{sec:analytical}, which is computed initially and after each accident event, two strategies fundamentally exist: vehicles can either gain access to network A from B starting at a critical time~$t^*$, or access is granted until time~$t^*$. Since the first variant would lead to unbounded total travel times for our modeled constant inflow and, consequently, to a diverging value of the objective function, only the second option qualifies as the optimal bang-bang control. Algorithmically, the computation of optimal control thus reduces to determining the optimal threshold~$t^*$, beyond which network access from B to A is prohibited.

To model the priority of accident avoidance, the weighting parameter $\lambda$ is introduced such that
\[
c_s = \lambda \cdot \int_0^T \left( f_A(t) + f_B(t) \right)\, dt.
\]
The parameter $\lambda$ thus represents the additional travel time per vehicle that the system operator is willing to accept to prevent one additional accident. For the weighting parameters of the objective function, we set $\theta = 0$ and analyze the values $\lambda \in \{0, \tfrac{1}{3}, \tfrac{2}{3}\}$, which corresponds to $c_s \in \{0, 4000, 8000\}$. This parameterization covers the range derived from established ececonimc appraisal figures for economic accident damages, the proportion of serious personal injuries, the value of travel time, and the value of statistical life.

\begin{figure}[H]
    \centering
    \includegraphics[width=0.8\textwidth]{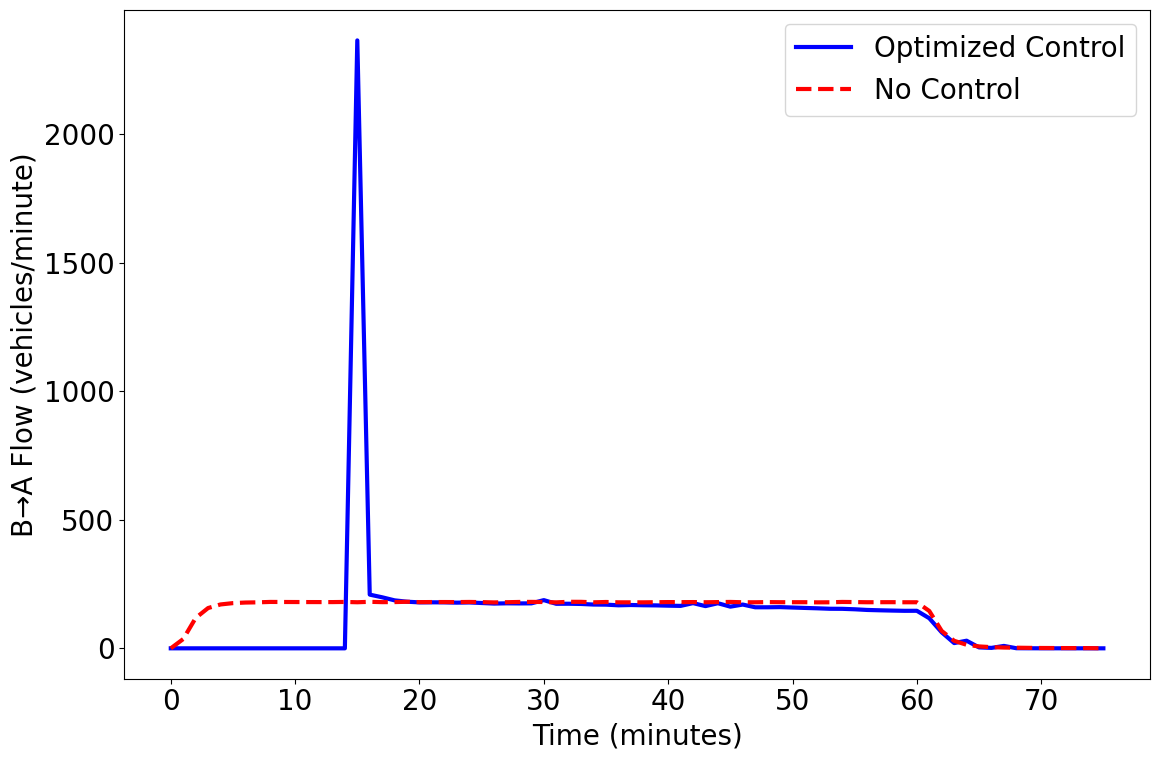}
    \caption{Transit flow from network B to A at an accident priority weight of $c_s = 1/3$ for the optimally controlled scenario and the uncontrolled reference scenario. After approximately 15 minutes, the optimally controlled system executes an abrupt release of all queued vehicles into network A, afterwards the flow stabilizes at a level slightly below the value observed in the uncontrolled reference scenario.}
    \label{fig:flow_comparison}
\end{figure}

The evaluation is based on 300 Monte Carlo simulations of the system. The following performance metrics were determined as arithmetic means across all realizations: Table~\ref{tab:travel_times} documents the average travel time per vehicle as well as the mean value of the objective function per vehicle, each with percentage changes relative to a no-control reference scenario. Table~\ref{tab:accidents} presents the average total number of accidents and the relative deviation from the uncontrolled scenario, also based on 300 simulation runs. Figure~\ref{fig:flow_comparison} visualizes the mean transit flow from network B to A for the case $c_s = 1/3$, both for the optimally controlled system and the no-control reference scenario.

\begin{table}[H]
\centering
\caption{Travel Time and Objective Function Results}
\label{tab:travel_times}
\renewcommand{\arraystretch}{1.3}
\begin{tabular}{lccc}
\toprule
\textbf{Scenario} & \textbf{Weight = 0.000} & \textbf{Weight = 0.333} & \textbf{Weight = 0.667} \\
\midrule
Base Accident Rate 
  & \makecell{(7.90, 7.90) \\ (-27.9\%, -27.9\%)} 
  & \makecell{(8.1, 9.1) \\ (-19.8\%, -20.4\%)} 
  & \makecell{(8.2, 10.2) \\ (-18.8\%, -20.1\%)} \\
High Accident Rate 
  & \makecell{(21.80, 21.80) \\ (-46.4\%, -46.4\%)} 
  & \makecell{(22.3, 25.3) \\ (-30.1\%, -29.7\%)} 
  & \makecell{(22.4, 28.0) \\ (-29.8\%,-30.1\%)} \\
\bottomrule
\end{tabular}
\end{table}

\begin{table}[h!]
\centering
\caption{Accident Results}
\label{tab:accidents}
\renewcommand{\arraystretch}{1.3}
\begin{tabular}{lccc}
\toprule
\textbf{Scenario} & \textbf{Weight = 0.000} & \textbf{Weight = 0.333} & \textbf{Weight = 0.667} \\
\midrule
Base Accident Rate 
  & 4.00, -13.1\% 
  & 3.06, -32.4\% 
  & 2.92, -35.4\% \\
High Accident Rate 
  & 11.00, -11.6\% 
  & 8.84, -28.0\% 
  & 8.61, -29.9\% \\
\bottomrule
\end{tabular}
\end{table}

The numerical results reveal the effectiveness of the proposed algorithm in managing the trade-off between travel time efficiency and safety across different accident scenarios and weighting parameters. 

The measured travel times are significantly reduced through implementation of the control algorithm in all scenarios, even when the control focus shifts more heavily towards accident prevention. The reduction in accident numbers ranges from 10\% to 35\% depending on the scenario, with the marginal additional benefit decreasing at higher $\lambda$ values. Remarkably, even at $\lambda = 0$, where accident reductions are not explicitly incorporated in the objective function, significant decreases in accident numbers are achieved. This results from the fact that accidents indirectly lead to extended travel times through degradation of the flow-density relationship. Hence, their number is minimized by the controller even when the objective function exclusively targets travel time minimization. Figure~\ref{fig:flow_comparison} illustrates that the optimal control exploits the nonlinear relationship between vehicle density and accident probability by initially withholding inflow to network A for approximately 15 minutes. This generates a pronounced peak in throughput at around 15 minutes. Subsequently, the optimally controlled throughput converges toward the uncontrolled reference value and deviates only moderately from it for the remaining simulation interval.

\section{Conclusion}
\label{sec:conclusion}

This paper presented a risk-aware perimeter control framework that explicitly couples travel time efficiency and safety objectives for large-scale urban networks. By integrating a self-exciting Hawkes process for accident modeling with the Generalized Bathtub Model dynamics, we developed a unified approach that captures both the stochastic nature of accident occurrences and their operational impacts on traffic flow.
The theoretical analysis revealed that optimal gating controls follow a bang-bang pattern with switches occurring only at accident events, which motivated an event-triggered MPC structure. This insight enabled significant computational savings—reducing optimization frequency from 60 times per hour to approximately 6 times per hour in the Copenhagen case study—while maintaining mathematical optimality. The risk-aware objective function, augmented with a variance penalty term, provides transportation agencies with a transparent mechanism to balance average performance against reliability. Our closed-form steady-state solutions demonstrate how risk aversion systematically biases the system toward safer operating configurations, with the single parameter $c_s$ offering intuitive control over this trade-off.
Numerical experiments on a network calibrated after metropolitan Copenhagen demonstrated the practical value of this integrated approach, achieving delay reductions up to 30\% and accident reductions up to 35\% relative to uncontrolled operations. Notably, even risk-neutral policies achieved significant safety improvements due to the explicit modeling of accident-induced capacity degradation, highlighting that safety and efficiency are inherently coupled in congested networks.
Several promising research directions emerge from this work. The two-reservoir abstraction could be extended to multi-region networks with more complex interaction patterns, though this would require careful analysis of the resulting switching logic and computational complexity. At a more granular level, the ideas behind the safety-aware control principles developed here could be integrated into the operational control of connected and automated vehicles, where the accident risk model could inform car-following behaviors and lane-changing decisions in real-time. This would create a multi-scale control architecture linking network-level perimeter control with vehicle-level trajectory planning.
The framework could also be enhanced by incorporating real-time accident prediction models using machine learning techniques to improve the anticipatory capabilities of the controller. From a methodological perspective, developing adaptive calibration methods for the Hawkes process parameters using streaming accident data would improve the controller's responsiveness to changing conditions. Finally, extending the approach to consider more comprehensive objectives—including environmental impacts, accessibility for different user groups, and economic costs beyond travel time—would create a more holistic optimization framework for urban mobility management. The modular structure of our framework, with its clear separation between dynamics, safety modeling, and optimization, provides a solid foundation for these future enhancements.

\newpage

\bibliographystyle{trb}
\bibliography{trb_template}

\end{document}